\documentclass[a4paper,12pt]{article}
\usepackage{amssymb}
\usepackage{latexsym}
\usepackage{amsmath,amsthm,amssymb}
\usepackage{amsfonts}

\newtheorem{tw}{Theorem}[section]

\newtheorem{cor}[tw]{Corollary}

\theoremstyle{definition}

\newtheorem{exa}[tw]{Example}
\newtheorem{rem}[tw]{Remark}

\begin{document}

\begin{center}
{\Large Penney's game between many players}
\end{center}
\begin{center}
{\sc Krzysztof Zajkowski}
\footnote{The author is supported by the Polish National Science Center, Grant no. DEC-2011/01/B/ST1/03838.}\\
Institute of Mathematics, University of Bialystok \\ 
Akademicka 2, 15-267 Bialystok, Poland \\ 
kryza@math.uwb.edu.pl 
\end{center}

\begin{abstract}
We recall a combinatorial derivation of the functions generating probability of winnings for each of many participants of the Penney's game and show a generalization of the Conway's formula to this case.
\end{abstract}

{\it 2010 Mathematics Subject Classification:} 61A06, 91A15 

\section{Introduction}
Let us  toss an  'unfair' coin with probabilities $p$ for heads ($H$) and $q=1-p$ for tails ($T$) and wait for an appearance of some chosen string of heads and tails. What is the expected number of tosses until this string occurs?  

Let $A=a_1a_2...a_l$ be a given pattern (a string  of heads and tails) of the length $l$. 
By $P(A)$ we will denote a value $P(a_1)P(a_2)\cdot...\cdot P(a_l)$. More precisely $P(A)$ is the probability of a cylindric set with  some fixed coordinates: $a_1$, $a_2$,... and $a_l$, respectively. We flip a coin until we get $A$ as a run  in the sequence of our trials. So we define a stopping time of the process in the following form
\begin{equation*}
\tau_A=\min\Big\{n\ge l:\;\xi_1\xi_2\ldots\xi_n,\;\xi_i\in\{H,T\}\; {\rm and} \;\xi_{n-l+1}=a_1,\xi_{n-l+2}=a_2,...,\xi_n=a_l\Big\},
\end{equation*}
if this minimum exists and $\infty$ if not. Now a more precise formulation of our question is: what is the expected value of $\tau_A$?
An answer was first given by Solov'ev in \cite[(1966)]{Sol}.
In the paper presented we show some combinatorial solution (compare Graham et. all \cite[VIII.8.4]{Concrete}), introducing at the same time notations required and presenting a model reasoning.

Let $A_n$ denote the set of  sequences in which the pattern $A$ appears exactly in the $n$-th toss, i.e. $A_n=\{\tau_A=n\}$, and $p_n$ the probability of $A_n$; $p_n=P(A_n)$.
Let $B_n$ denote a set of sequences in which $A$ does not appear in the first $n$ tosses, i.e. $B_n=\{\tau_A>n\}$, and its probability by $q_n=P(B_n)$.
Let us consider now a set of sequences in which $A$ does not appear in the first $n$ tosses  and  appears in the next $l$ trials, i.e. the set
\begin{equation*}
\Big\{(\xi_k)\in\{H,T\}^\mathbb{N}:\;(\xi_k)\in B_n\;{\rm and}\; \xi_{n+1}=a_1,\xi_{n+2}=a_2,...,\xi_{n+l}=a_l\Big\}.  
\end{equation*}
It seems
that the probability of the set amounting to $q_nP(A)$ is equal to $p_{n+l}$  but we must check whether $A$ does not occur earlier in the trials from $n+1$ to $n+l-1$.

Let $A_{(k)}$ and $A^{(k)}$  denote  strings of $k$-first and $k$-last terms of $A$ ($1\le k \le l$), respectively. Note that $A_{(l)}=A^{(l)}=A$. 
Let $[A_{(k)}=A^{(k)}]$ equals $1$ if $A_{(k)}=A^{(k)}$ or $0$ if not. Additionally, let us assume that $P(A^{(0)})=1$. Now we can write the formula on
$q_nP(A)$ as follows:
\begin{equation}
\label{equ1}
q_nP(A)=\sum_{k=1}^{l}[A_{(k)}=A^{(k)}]P(A^{(l-k)})p_{n+k}.
\end{equation}
Observe that the $l$-th summand in the above is equal to $p_{n+l}$. Remembering that $p_0=p_1=...=p_{l-1}=0$, multiplying the above equation by $s^{n+l}$ and summing from $n=0$ to infinity we get
\begin{equation}
\label{equ2}
Q_{\tau_A}(s)P(A)s^l=g_{\tau_A}(s)\sum_{k=1}^{l}[A_{(k)}=A^{(k)}]P(A^{(l-k)})s^{l-k},
\end{equation}
where $g_{\tau_A}(s)=\sum_{n=0}^\infty p_ns^n$ is the probability generating function for a random variable $\tau_A$ of number of tosses until $A$ occurs and $Q_{\tau_A}(s)=\sum_{n=0}^\infty q_ns^n$ is the generating function of tail probabilities $q_n$. Because one can bound $\tau_A$ by a random variable with  the geometric distribution then one can show that $E\tau<\infty$. Hence $q_n=\sum_{k=n+1}^\infty p_k$ so we can obtain the second equation that relates $g_{\tau_A}$ and $Q_{\tau_A}$:
\begin{equation}
\label{equ3}
Q_{\tau_A}(s)=\frac{1-g_{\tau_A}(s)}{1-s}.
\end{equation}
Solving the above two equations we get
$$
g_{\tau_A}(s)=\frac{P(A)s^l}{P(A)s^l+(1-s)\sum_{k=1}^{l}[A_{(k)}=A^{(k)}]P(A^{(l-k)})s^{l-k}}
$$
and
$$
Q_{\tau_A}(s)=\frac{\sum_{k=1}^{l}[A_{(k)}=A^{(k)}]Pr(A^{(l-k)})s^{l-k}}{P(A)s^l+(1-s)\sum_{k=1}^{l}[A_{(k)}=A^{(k)}]P(A^{(l-k)})s^{l-k}}. 
$$
Since $E\tau_A=Q_{\tau_A}(1)$, one can calculate the general formula for the expected number of tosses as follows:
$$
E\tau_A=\frac{\sum_{k=1}^{l}[A_{(k)}=A^{(k)}]P(A^{(l-k)})}{P(A)}=\sum_{k=1}^{l}\frac{[A_{(k)}=A^{(k)}]}{P(A_{(k)})}.
$$
It is the answer to the question posed at the beginning.
\begin{rem}
Note that the above solution can be also obtained by using  Renewal Theorem (see \cite[Lem.1]{Chen}) and a martingale approach, called the gambling team method, (see \cite[Lem.2.4]{Li}).
\end{rem}
In the classical Penney Ante game (see \cite{Pen}), for a given string of  fix length we want to show a second one of the same length with a higher probability to be the first to occur. In \cite{Chen}, Chen and Zame proved that for two-person games, public knowledge of the opponent's string leads to an advantage.
Guibas and Odlyzko \cite{Gui} show some optimal strategy for second player.
An algorithm for computing the odds of winning for the competing patterns was discovered by Conway and described by Gardner \cite{Con}.
The Conway's formula  allows us to compare  probability of winnings for  two players.

The main aim of this paper is to show a generalization of the Conway's formula to the case of many gamblers (Section 3).
But first (Section 2) we present a derivation and solution of  the system of equations proposed by  Guibas and Odlyzko  \cite[Th.3.3]{Gui}, and in a similar and much general form (with some starting string) by Li and Gerber \cite[the system (39)]{GerLi}.

\section{Functions  generating probability of winnings}
Let $m$ players choose $m$ strings $A_i$ ($1\le i \le m$) of heads and tails of lengths $l_i$, respectively. We start to toss an 'unfair' coin and
wait for an occurrence of some $A_i$.
We ask about the chances of  winning for each player; that is about the probability $p_{A_i}$ that the string $A_i$ to be the first to occur.
We assume that any $A_i$ is not a substring of other $A_j$; in an opposite case $p_{A_j}=0$ or for some sequences both players may win simultaneously.

Let $\tau$ denote the number of tosses to the end of the game, i.e. $\tau=\min\{\tau_{A_i}:\;1\le i \le m\}$, where
$\tau_{A_i}$ is the stopping time until pattern $A_i$ occurs.  Notice that $P(\tau=n)=p_n=\sum_{i=1}^m p_n^{A_i}$, where $p_n^{A_i}=P(\tau=\tau_{A_i}=n)$ is the probability
that the $i$-th player wins exactly in the $n$-th toss. Let $g_\tau$ and $g_\tau^{A_i}$ denote the functions  generating  distributions of probability
$(p_n)$ and $(p^{A_i}_n)$, respectively, and $Q_\tau$ the generating function of tail distributions  $q_n=P(\tau>n)$.

Similarly, as in Introduction, let $B_n=\{\tau>n\}$ be the set of sequences of tails and heads in which any string $A_i$ does not appear in first $n$ tosses.
In the system of $m$ patterns if we add the string $A_i$ to the set $B_n$  then we must check
if neither $A_i$ nor other patterns appear earlier. 
For this reason a system of equations   
$$
q_nP(A_i)=\sum_{j=1}^m\sum_{k=1}^{\min\{l_i,l_j\}}[A_{i(k)}=A_j^{(k)}]P(A_i^{(l_i-k)})p^{A_j}_{n+k},
$$
for each $1\le i \le m$, where $[A_{i(k)}=A_j^{(k)}]=1$ if $A_{i(k)}=A_j^{(k)}$ or $0$ if not,
corresponds to the equation (\ref{equ1}).

Multiplying the above equation by $s^{n+l_i}$ and summing from $n=0$ to infinity we get the following recurrence equations
\begin{equation*}
Q_\tau(s)P(A_i)s^{l_i}  =  \sum_{j=1}^m g_\tau^{A_j}(s)\sum_{k=1}^{\min\{l_i,l_j\}}[A_{i(k)}  =  A_j^{(k)}]P(A_i^{(l_i-k)})s^{l_i-k}.
\end{equation*}
Let $w_{A_i}^{A_j}(s)$ denote the polynomial $\sum_{k=1}^{\min\{l_i,l_j\}}[A_{i(k)}  =  A_j^{(k)}]P(A_i^{(l_i-k)})s^{l_i-k}$; now we can rewrite the
above system of $m$ equations as follows
\begin{equation}
\label{sys1}
Q_\tau(s)P(A_i)s^{l_i}  =  \sum_{j=1}^m g_\tau^{A_j}(s)w_{A_i}^{A_j}(s)\quad (1\le i \le m).
\end{equation}
Since $g_\tau =\sum_{j=1}^m g_\tau^{A_j}$, by virtue of (\ref{equ3}), we get
$$
Q_\tau(s)=\frac{1-\sum_{j=1}^m g_\tau^{A_j}(s)}{1-s};
$$
compare the above equation and (\ref{sys1}) with the system of linear equations given by Guibas and Odlyzko in Theorem 3.3  \cite{Gui}.

Inserting the form of $Q_\tau$ into (\ref{sys1}) we obtain
$$
P(A_i)s^{l_i}=\sum_{j=1}^m g_\tau^{A_j}(s)[P(A_i)s^{l_i}+(1-s)w_{A_i}^{A_j}(s)]\quad (1\le i \le m).
$$
To solve this system of functional equations we use the Cramer's rule. 
Define now a functional matrices  

$$
\mathcal{A}(s) =
\begin{matrix}
P(A_i)s^{l_i}+(1-s)w_{A_i}^{A_j}(s)
\end{matrix}
_{1\le i,j \le m}
$$ 
and
\begin{equation}
\label{matrB}
\mathcal{B}(s)=
\begin{pmatrix}
w_{A_i}^{A_j}(s)
\end{pmatrix}
_{1\le i,j \le m}.
\end{equation}
Notice that because $w_{A_i}^{A_i}(0)=1$, and $w_{A_i}^{A_j}(0)=0$ for $i\neq j$ then $\mathcal{A}(0)$ and $\mathcal{B}(0)$ are the identity matrices.
Let $\mathcal{B}^j(s)$ denote the matrix formed by replacing the $j$-th column of $\mathcal{B}(s)$ by  the column vector $[P(A_i)s^{l_i}]_{1\le i \le m}$. Because the determinant 
of matrices $m\times m$ is a $m$-linear functional with respect to columns (equivalently to rows) then one can check that 
$$
\det \mathcal{A}(s)=(1-s)^m\det \mathcal{B}(s) +(1-s)^{m-1}\sum_{j=1}^m \det \mathcal{B}^j(s).
$$
The determinant $\det\mathcal{A}(s)$ is a polynomial of variable $s$ and $\det\mathcal{A}(0)=1$. For these reasons $\det\mathcal{A}(s)\neq 0$ in some neighborhood of zero.
It means that in this neighborhood there exists a solution of the system.

If now, similarly, $\mathcal{A}^j(s)$ denotes the matrix formed by replacing the $j$-th column of $\mathcal{A}(s)$ by  the column vector $[P(A_i)s^{l_i}]_{1\le i \le m}$, then the determinant's
calculus gives $\det \mathcal{A}^j(s)=(1-s)^{m-1}\det \mathcal{B}^j(s)$.  Finally, by the  Cramer's rule we obtain:
$$
g_\tau^{A_i}(s)=\frac{\det\mathcal{A}^i(s)}{\det\mathcal{A}(s)}=\frac{\det \mathcal{B}^i(s)}{\sum_{j=1}^m \det \mathcal{B}^j(s)+(1-s)\det \mathcal{B}(s) }
$$
for $1\le i \le m$.
In this way we have proved the following 
\begin{tw}
\label{Mtw}
If $m$ players choose $m$ strings of heads and tails $A_i$ $(1\le i \le m)$ such that any $A_i$ is not a substring of other $A_j$ then
the function $g_\tau^{A_i}$ generating probability   of winning of the $i$-th player is given by the following formula:
\begin{equation}
\label{pgf}
g_\tau^{A_i}(s)=\frac{\det \mathcal{B}^i(s)}{\sum_{j=1}^m \det \mathcal{B}^j(s)+(1-s)\det \mathcal{B}(s)},
\end{equation}
where $\mathcal{B}(s)$ is the matrix defined by (\ref{matrB}).
\end{tw}

Notice that the probability generating function $g_\tau^{A_i}(s)=\sum_{n=0}^\infty p_n^{A_i}s^n$ is undoubtedly well define on the interval $[-1,1]$
(it is an analytic function on $(-1,1)$). The right hand side of (\ref{pgf}) is a rational function equal to $g_\tau^{A_i}$ in some neighborhood of zero.
By the analytic extension we know that there exists the limit of the right hand side of (\ref{pgf}) by $s\to 1^-$ which is equal to $g_\tau^{A_i}(1)$. 
Thus the probability $p_{A_i}$ that the string $A_i$ occurs  first is given by the following formula 
\begin{equation}
\label{pro}
g_\tau^{A_i}(1)=\frac{\det \mathcal{B}^i(1)}{\sum_{j=1}^m \det \mathcal{B}^j(1)},
\end{equation}
where the right hand side of the above equation is understood as the limit of (\ref{pgf}) with $s$ approaching  the left-side to 1 ($s\to 1^-$).

\section{A generalization of the Conway's formula}

Define a number $A_j\ast A_i$ as 
$$
\frac{w_{A_i}^{A_j}(1)}{P(A_i)}=\sum_{k=1}^{\min\{l_i,l_j\}}[A_{i(k)}=A_j^{(k)}]\frac{P(A_i^{(l_i-k)})}{P(A_i)}
=\sum_{k=1}^{\min\{l_i,l_j\}}\frac{[A_{i(k)}=A_j^{(k)}]}{P(A_{i(k)})};
$$
compare \cite[Def.2.2]{Li}.
Define now a matrix 
$
\mathcal{C}=
\begin{pmatrix}
(A_j\ast A_i)
\end{pmatrix}
_{1\le i,j \le m}.
$
Observe that $\det\mathcal{B}(1)=\prod_{i=1}^m P(A_i)\det\mathcal{C}$ and $\det\mathcal{B}^j(1)=\prod_{i=1}^m P(A_i)\det\mathcal{C}^j$, where $\mathcal{C}^j$
is the matrix formed by replacing the $j$-th column of $\mathcal{C}$ by  the column vector of units
$\begin{matrix} 1\end{matrix}_{1\le i \le m}$. Due to (\ref{pro}) and the above observations we can formulate the following
\begin{cor}
\label{Mcor}
The probability that the $i$-th player wins is equal to
$$
p_{A_i}=\frac{\det \mathcal{C}^i}{\sum_{j=1}^m \det \mathcal{C}^j}.
$$
\end{cor}
Let us emphasize that the above Corollary is the generalization of the Conway's formula. For two players we get
$$
\frac{p_{A_1}}{p_{A_2}}=\frac{\det\mathcal{C}^1}{\det\mathcal{C}^2}=
\det
\begin{pmatrix}
1 & (A_2\ast A_1)\\
1 & (A_2\ast A_2)
\end{pmatrix}
:
\det
\begin{pmatrix}
(A_1\ast A_1) & 1\\
(A_1\ast A_2) & 1
\end{pmatrix}
=\frac{(A_2\ast A_2)-(A_2\ast A_1)}{(A_1\ast A_1)-(A_1\ast A_2)}.
$$

On  page 188 in \cite{Gui} one can find a expanded form of the probability of winning of a player $A$ in  three-person $A,B,C$ game.
In the following example we present an applications of our results to such case.
\begin{exa}
Take three strings of heads  and tails : $A_1=THH$, $A_2=HTH$ and $A_3=HHT$. In this case 
\begin{equation*}
\mathcal{B}(s)=
\begin{pmatrix}
w_{A_i}^{A_j}(s)
\end{pmatrix}
_{1\le i,j \le 3}
=
\begin{pmatrix}
1 & ps & p^2s^2 \\
pqs^2 & pqs^2+1 & ps\\
pqs^2+qs & pqs^2 & 1
\end{pmatrix}.
\end{equation*}
By Theorem \ref{Mtw} one can obtain the probability generating functions for winnings of $i$-th player. The matrix 
\begin{equation*}
\mathcal{C}=
\begin{pmatrix}
\frac{w_{A_i}^{A_j}(1)}{P(A_i)}
\end{pmatrix}
_{1\le i,j \le 3}
=
\begin{pmatrix}
\frac{1}{p^2q} & \frac{1}{pq} & \frac{1}{q} \\
\frac{1}{p} & \frac{pq+1}{p^2q} & \frac{1}{pq}\\
\frac{p+1}{p^2} & \frac{1}{p} & \frac{1}{p^2q}
\end{pmatrix}.
\end{equation*}
and 
\begin{equation*}
\mathcal{C}^1=
\begin{pmatrix}
1 & \frac{1}{pq} & \frac{1}{q} \\
1 & \frac{pq+1}{p^2q} & \frac{1}{pq}\\
1 & \frac{1}{p} & \frac{1}{p^2q}
\end{pmatrix},\quad
\mathcal{C}^2=
\begin{pmatrix}
\frac{1}{p^2q} & 1 & \frac{1}{q} \\
\frac{1}{p} & 1 & \frac{1}{pq}\\
\frac{p+1}{p^2} & 1 & \frac{1}{p^2q}
\end{pmatrix},\quad
\mathcal{C}^3=
\begin{pmatrix}
\frac{1}{p^2q} & \frac{1}{pq} & 1\\
\frac{1}{p} & \frac{pq+1}{p^2q} & 1\\
\frac{p+1}{p^2} & \frac{1}{p} & 1
\end{pmatrix}.
\end{equation*}
On the basis of Corollary \ref{Mcor} we can calculate the probability that the $i$-th player wins:
$$
p_{A_1}=\frac{q(1+pq)}{1+q},\quad p_{A_2}=\frac{q}{1+q},\quad p_{A_3}=\frac{p(1-q^2)}{1+q}.
$$
For the fair coin $p_{A_1}=\frac{5}{12}$, $p_{A_2}=\frac{1}{3}$ and $p_{A_3}=\frac{1}{4}$.
\end{exa}

\begin{rem}
It is possible to consider a situation when we get to know of the first outcomes of the game. It leads to considerations 
of the game with an initial pattern (compare Li \cite[the system of equations (3.7)]{Li}). A form of solutions of that system  presented in \cite[Cor.4.3]{Zaj} may also be  treated as some another generalization of the Conway's formula to the case when an initial pattern is known.
\end{rem}



\begin{thebibliography}{                    }

\bibitem{Chen}
Chen R. and Zame A. (1979) {\em On the fair coin-tossing games}, J. Multivariate Anal., vol. 9, 150-157. 

\bibitem{Con}
Gardner M. (1974) {\em On the paradoxical situations that arise from nontransitive relations}, Scientific American 231 (4), 120-124.

\bibitem{GerLi}
Gerber H.\,U. and Li S-Y.\,R. (1981) {\em The occurrence of sequence patterns in repeated experiments and hitting times in a Markov chain}, Stochastic Processes and Their Applications 11, 101-108.

\bibitem{Gui}
Guibas L.\,J. and Odlyzko A.\,M.(1981) {\em String overlaps, pattern matching, and nontransitive games}, Journal of Combinatorial Theory (A) 30, 183-208.


\bibitem{Concrete}
Graham R.\,L.,  Knuth D.\,E. and  Patashnik O. (1989) {\em Concrete Mathematics: a Foundation for Computer Science}, Addison-Wesley Publishing Company.

\bibitem{Li}
Li S-Y.\,R. (1980) {\em A martingale approach to the study of occurrence of sequence patterns in repeated experiments}, The Annals of Probability, Vol. 8, 1171-1176. 

\bibitem{Pen} 
Penney W. (1974) {\em Problem 95: Penney-Ante}, Journal of Recreational Mathematics 7, 321. 

\bibitem{Sol}
Solov'ev A.\,D.(1966) {\em A combinatorial identity and its application to the problem concerning the first occurrence of a rare event},
Theory of Probability and its Applications 11, 313-320.

\bibitem{Zaj}
Zajkowski K. (2014) {\em A note on the gambling team method}, Statistics and Probability Letters 85, 45-50.

\end{thebibliography}
\end{document}